\documentclass{mva_style}
\usepackage[T1]{fontenc}
\usepackage{graphicx}
\usepackage[cmex10]{amsmath}
\usepackage{color}
\usepackage{amsmath,amsfonts}
\usepackage{subcaption}
\usepackage{epstopdf}
\usepackage{bbm}
\usepackage{float}
\usepackage[ruled,vlined]{algorithm2e}
\DeclareMathOperator*{\argmin}{argmin}
\newtheorem{theorem}{Theorem}

\newcommand{\myindent}[1]{
	\newline\makebox[#1cm]{}
}

\finalcopy 

\begin{document}
\title{A Fast Algorithm for a Weighted Low Rank Approximation}

\author{Aritra Dutta\\
Department of Mathematics\\
University of Central Florida\\
Orlando, FL 32816\\
{\tt d.aritra2010@knights.ucf.edu}\\
\and
Xin Li\\
Department of Mathematics\\
University of Central Florida\\
Orlando, FL 32816\\
{\tt xin.li@ucf.edu}\\
}

\maketitle

\section*{\centering Abstract}
\vspace{-0.01in}
\textit{
   Matrix low rank approximation including the classical PCA and the robust PCA (RPCA) method have been applied to solve the background modeling problem in video analysis.~Recently, it has been demonstrated that a special weighted low rank approximation of matrices can be made robust to the outliers similar to the $\ell_1$-norm in RPCA method.~In this work, we propose a new algorithm that can speed up the existing algorithm for solving the special weighted low rank approximation and demonstrate its use in background estimation problem.
}
\vspace{-0.05in}
\section{Introduction}
\vspace{-0.05in}
Background estimation is one of the crucial steps in video analysis systems.~The celebrated eigen-background model proposed in \cite{oliver1999} was the first case when principal component analysis~(PCA) was used in background modeling.
Recently, using a sparse model for the foreground and a low rank model for the background, \cite{candeslimawright,APG} proposed a robust principal component analysis~(RPCA) model.

For an integer $r\le\min\{m,n\}$ and a matrix $A\in\mathbb{R}^{m\times n}$, the classical PCA can be cast as:
\vspace{-0.09in}
\begin{eqnarray}
\label{pca}
\min_{\substack{{X}\in\mathbb{R}^{m\times n}\\{\rm r}({X})\le r}}\|A-{X}\|_F^2,
\end{eqnarray}
~\\[-0.2in]
where ${\rm r}({X})$ denotes the rank of the matrix ${X}$ and $\|\cdot\|_F$ denotes the Frobenius norm of matrices. The solutions to~(\ref{pca}) are given using the singular value decompositions~(SVDs) of $A$ through the hard thresholding operations on the singular values:
\vspace{-0.1in}
\begin{align}\label{hardthresholding}
X^*=H_r:=U\Sigma_rV^T,
\end{align}
~\\[-0.24in]
where
$
A=U\Sigma V^T,
$
is a SVD of $A$ and $\Sigma_r$ is the diagonal matrix obtained from $\Sigma$ by hard-thresholding:~keeping only the largest $r$ entries and replacing the others by 0.~In literature this is also referred to as Eckart-Young-Mirsky's theorem~\cite{pca}.~The solutions to~(\ref{pca}) as given in~(\ref{hardthresholding}) suffer from the fact that none of the entries of $A$ is guaranteed to be preserved in $X^*$~\cite{lupeiwang,shpak}.~This could be a limitation of PCA, as in many real world problems one has good reasons to keep certain entries of $X$ unchanged while looking for a low rank approximation.~For example, if we know that certain frames of the input video matrix $A$, say frames \#1 and \#5, are pure background, then we may insist on preserving columns \#1 and \#5 when looking for a low rank approximation.~In 1987, Golub, Hoffman, and Stewart proposed the following {\it constrained} low rank approximation problem~\cite{golub}: Given $A=(A_1\;\;A_2)\in\mathbb{R}^{m\times n},$ find $\tilde{A}_2$ such that
~\\[-0.23in]
\begin{eqnarray}
\tilde{A}_2=\argmin_{X_2:{\rm r}(A_1\;\;X_2)\le r}\|(A_1\;\;A_2)-(A_1\;\;{X}_2)\|_F^2.\label{golub's problem}
\end{eqnarray}
That is, Golub, Hoffman, and Stewart required that the first few columns, $A_1,$ of $A$ must be preserved when one looks for a low rank approximation of $(A_1\;\;A_2).$  As in the standard low rank approximation, the constrained low-rank approximation problem of Golub, Hoffman, and Stewart also has a closed form solution.
~\\[-0.1in]
\begin{theorem}\cite{golub}
\label{theorem 1}
If $A = (A_1\;\;A_2)\in\mathbb{R}^{m\times n}$ with $k={\rm r}(A_1)$ and $r\ge k,$ then the solutions to~(\ref{golub's problem}) are given by
\vspace{-0.1in}
	\begin{align}\label{ghs}
	\tilde{A_2}= P_{A_1}(A_2)+H_{r-k}\left(P^{\perp}_{A_1}(A_2)\right),
	\end{align}
	~\\[-0.2in]
where $P_{A_1}$ and $P^\perp_{A_1}$ are the projection operators to the column space of $A_1$ and its orthogonal complement, respectively.
\end{theorem}

Instead of requiring exact matching in the first few columns,~as in~problem~(\ref{golub's problem}), we may only be interested in the case when the first few columns are close to the given ones.~For example,
for background estimation, we may have prior knowledge that some frames (say, represented by the columns of $A_1$) are {\it almost} pure background. So, 
\vspace{-0.09in}
$$
A_1=A_1^*+E$$
~\\[-0.24in]
for some true background frames $A_1^*$ and {\it small} noise $E$. Thus, we need to recover $A_1^*$ instead of matching $A_1$ exactly. So, we consider the following problem:
Given $A=(A_1\;\;A_2)\in\mathbb{R}^{m\times n}$ and $W_1\in\mathbb{R_+}^{m\times k}$, solve:
\vspace{-0.05in}
\begin{eqnarray}\label{2nd approximate golub's problem}
\min_{\substack{X_1,X_2\\{\rm r}(X_1\;\;X_2)\le r}}\{\|(A_1-X_1)\odot W_1\|^2_F+\|A_2-X_2\|_F^2\},
\end{eqnarray}
~\\[-0.2in]
where $\odot$ denotes the entrywise multiplication.~Problem~(\ref{2nd approximate golub's problem}) is a special case of weighted low rank approximation~\cite{srebro,manton,markovosky}:~Consider the following problem with $W=(W_1\;\;W_2)$ of non-negative terms
\vspace{-0.05in}
\begin{eqnarray}\label{hadamard problem}
	\min_{\substack{X_1,X_2\\{\rm r}(X_1\;\;X_2)\le r}}\|\left((A_1\;\;A_2)-({X}_1\;\;{X}_2)\right)\odot(W_1\;W_2)\|_F^2.
\end{eqnarray}
~\\[-0.15in]
Unlike classical (unweighted) low rank approximation, problem~(\ref{hadamard problem}) has no closed form solution in general~\cite{srebro}. So, numerical methods must be employed.~Recently, it has been demonstrated in~\cite{duttali_cvpr} that a method based on solving~(\ref{2nd approximate golub's problem}) can outperform the RPCA methods.~In this paper we propose a faster algorithm by exploiting an interesting property of the solution to problem~(\ref{2nd approximate golub's problem}).~Our algorithm is capable of achieving the desired accuracy faster as compared to~\cite{duttali,duttali_cvpr} and outperforming RPCA methods~\cite{LinChenMa,candeslimawright,APG} in background estimation problem.

The rest of the paper is organized as follows.~In Section~\ref{property},~we make an important observation on the solution to~(\ref{2nd approximate golub's problem}).~Based on this observation, we propose a new algorithm to solve problem~(\ref{2nd approximate golub's problem}) in Section~\ref{algorithm}.~Numerical results demonstrating the performance of the proposed algorithm are given in Section~\ref{numerical}.
\section{An Interesting Observation}\label{property}
\vspace{-0.05in}
We will design our algorithm using the observation as stated in the following result.
\begin{theorem}\label{theorem 7}
	Assume $r>k$.~For $(W_1)_{ij}>0$, if $(\hat{X}_1,\hat{X}_2)$ is a solution to~(\ref{2nd approximate golub's problem}), then
	\vspace{-0.1in}
	$$
	\hat{X}_2= P_{\hat{X}_1}(A_2)+H_{r-k}\left(P^{\perp}_{\hat{X}_1}(A_2)\right).
	$$
\end{theorem}
\vspace{-0.05in}
{\it Proof.} Note that,
\begin{align*}
&\|(A_1-\hat{X}_1)\odot W_1\|_F^2+\|A_2-\hat{X}_2\|_F^2\\
=&\min_{\substack{X_1,X_2\\{\rm r}(X_1\;\;X_2)\le r}}\left(\|(A_1-X_1)\odot W_1\|_F^2+\|A_2-X_2\|_F^2\right)\\
\le&\|(A_1-\hat{X}_1)\odot W_1\|_F^2+\|A_2-{X}_2\|_F^2,
\end{align*}
for all $X_2$ such that ${\rm r}(\hat{X}_1\;\;X_2)\le r.$
Thus,
\begin{align}
\hat{X}_2=\arg\min_{\substack{{X}_2\\{\rm r}(\hat{X}_1\;\;X_2)\le r}}\|A_2-X_2\|_F^2.
\end{align}
Thus, applying Theorem~\ref{theorem 1} with $A_1=\hat{X}_1$, we get
$\hat{X}_2= P_{\hat{X}_1}(A_2)+H_{r-k}\left(P^{\perp}_{\hat{X}_1}(A_2)\right).$ ~$\Box$

\vspace{-0.07in}
\section{Algorithm}\label{algorithm}
\begin{figure*}
 	\centering
 	\begin{subfigure}{.65\columnwidth}
 		\includegraphics[width=\columnwidth]{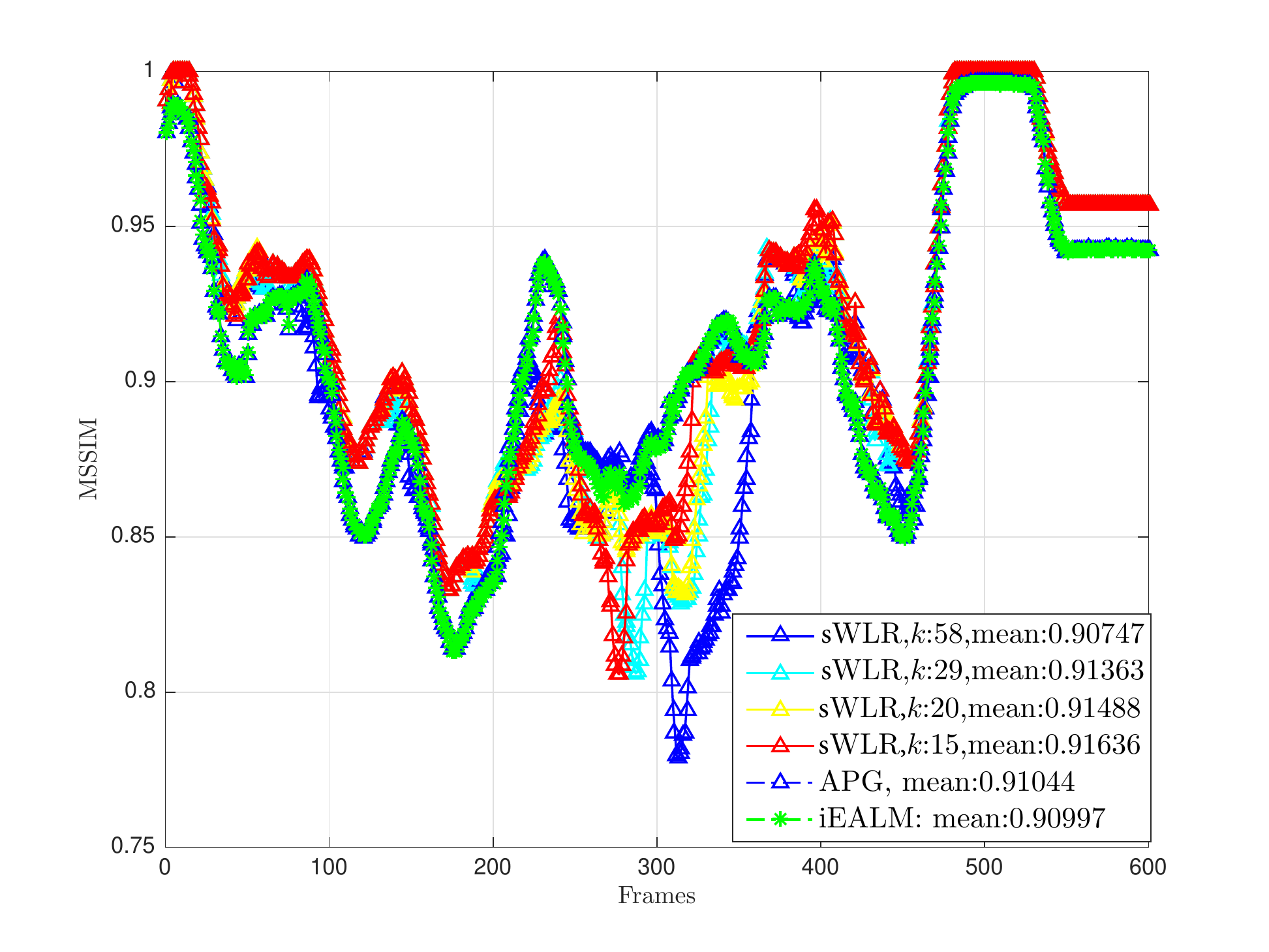}%
 		\caption{Basic}
 	\end{subfigure}\hfill
 	\begin{subfigure}{.65\columnwidth}
 		\includegraphics[width=\columnwidth]{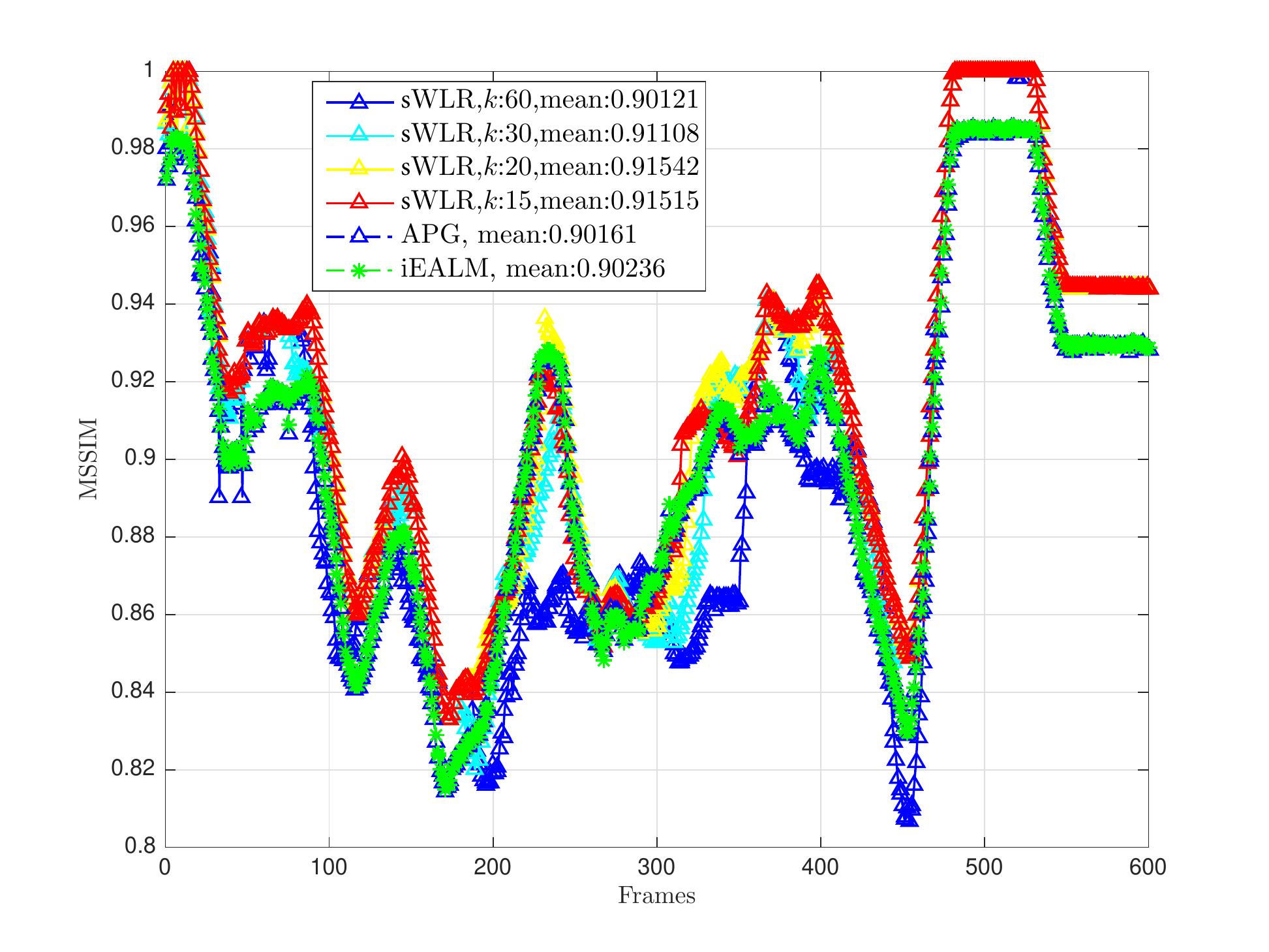}%
 		\caption{Noisy night}
 \end{subfigure}\hfill
 		\begin{subfigure}{.65\columnwidth}
 			\includegraphics[width=\columnwidth]{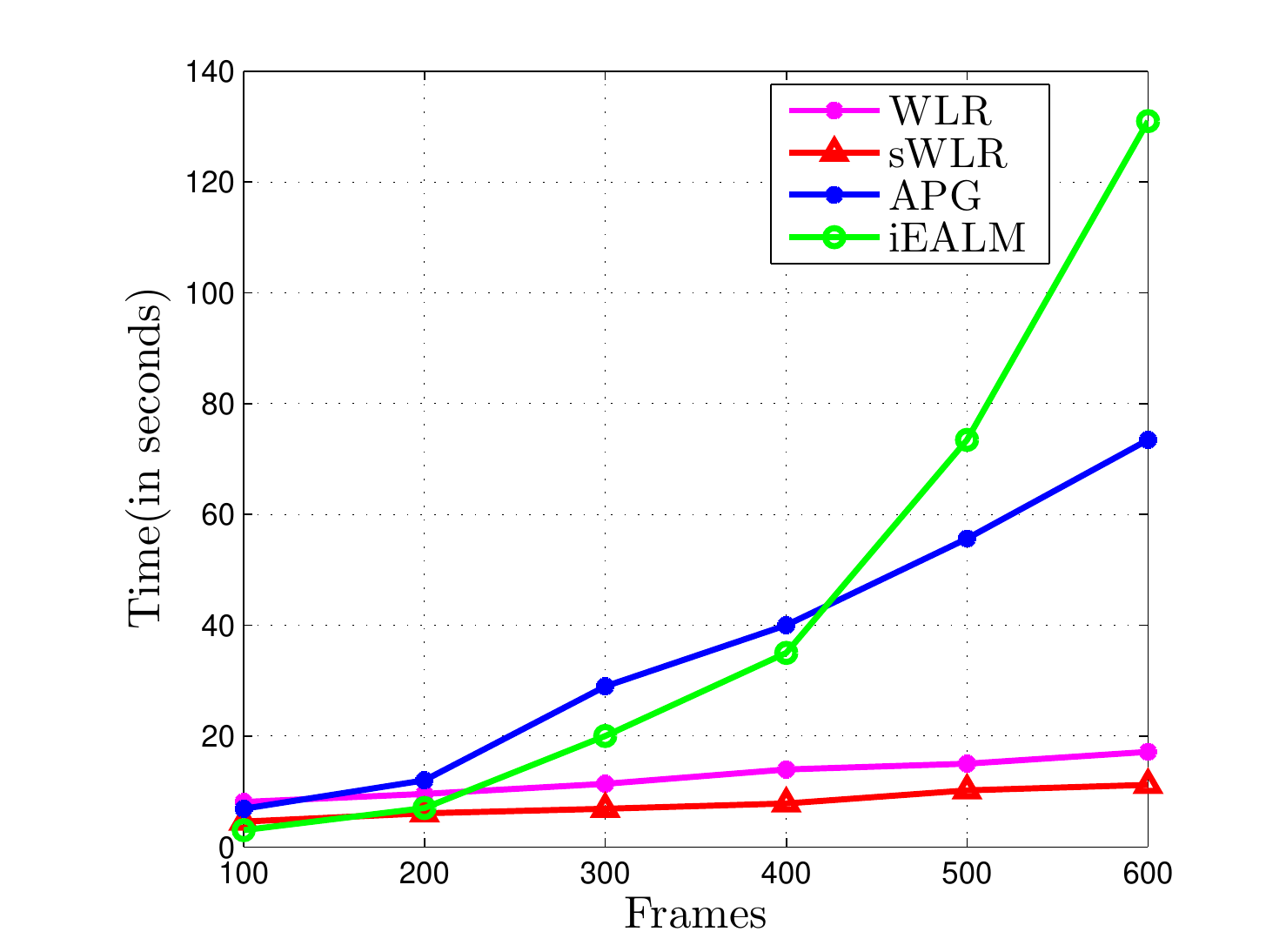}%
 			\caption{}
 			 		\end{subfigure}
 	\caption{\small Mean SSIM of different methods on: (a)~{\it Basic} and (b)~{\it Noisy night} scenario.~The choice of $k$ for sWLR is critical and empirically determined.~However, sWLR has better MSSIM compare to RPCA algorithms corresponding to the frames which has static foreground~(frame numbers 551 to 600) or no foreground~(frame numbers 6 to 12 and 483 to 528).~(c)~Number of video frames vs.~computation time on the {\it Basic} sequence.~As the number of video frame increases RPCA algorithms show an exponential increment in computational time.~For both WLR and sWLR, $k=15$, $r=k+1$.}
 	\label{ssim}
 \end{figure*}

\vspace{-0.09in}
In this section we propose a numerical algorithm to solve~(\ref{2nd approximate golub's problem}).~We do not use the general algorithm as in~\cite{srebro,wibergjapan,wiberg,markovosky} for solving~(\ref{hadamard problem}), since we focus on the special weight where $W_2=\mathbf{1}$, a matrix of all 1s.~In~\cite{duttali,duttali_cvpr},~the authors proposed an algorithm WLR to solve~(\ref{2nd approximate golub's problem}) which takes advantage of the special weight and performs much faster than the general weighted algorithm.~A rigorous comparison of accuracy and efficiency of WLR compare to the general weighted low rank approximation algorithms is discussed in~\cite{duttali}.~In this Section, we propose an accelerated version of the algorithm proposed in~\cite{duttali,duttali_cvpr}~(see Figure~\ref{ssim}~(c)) and demonstrate its use in background estimation.
\begin{figure}
	\begin{center}
		\includegraphics[width=\linewidth]{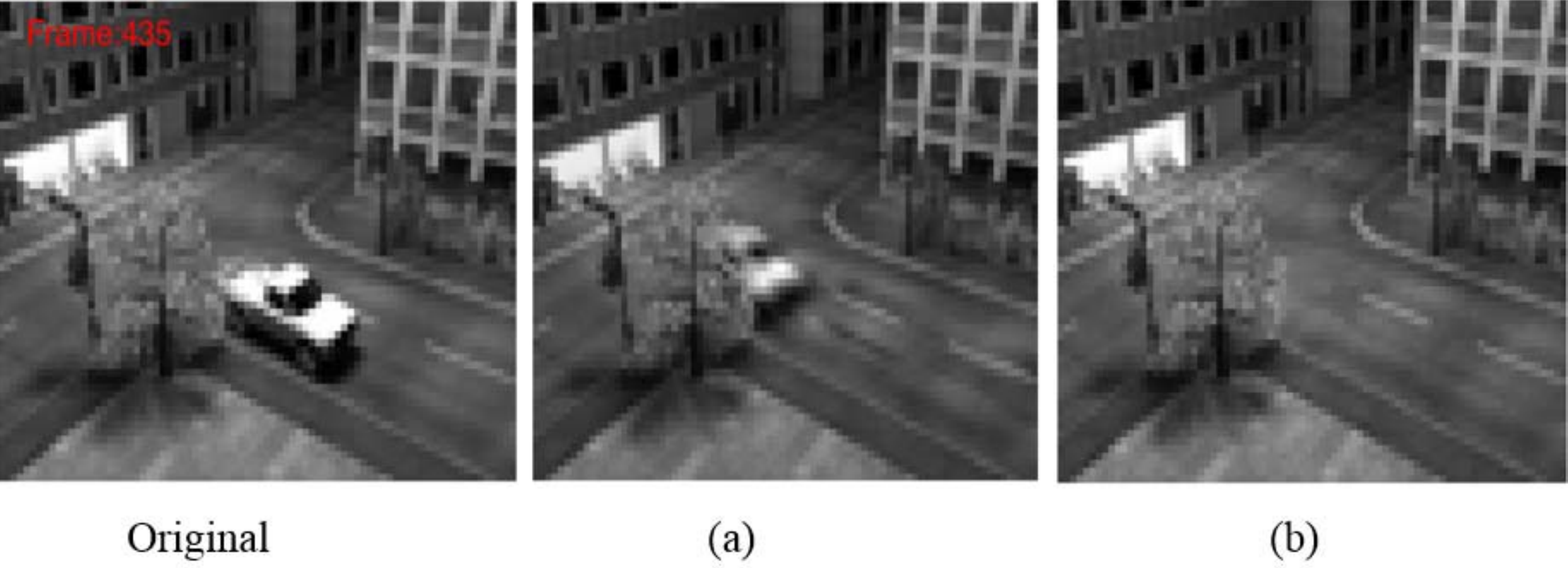}
	\end{center}
	\caption{The effect of using weights in sWLR algorithm on the~{\it Basic} scenario.~Frame number 435.~Background estimation using sWLR with:~(a)~$(W_1)_{ij}\in[5,10]$,~(b)~$(W_1)_{ij}\in[500,1000].$~In~(a)~the estimated background has blurry patches of the foreground object, but as we increase the weights, the foreground object disappears in~(b).}
	\label{weight_435}
\end{figure}
 \begin{figure}
 	\centering  \includegraphics[width=\linewidth]{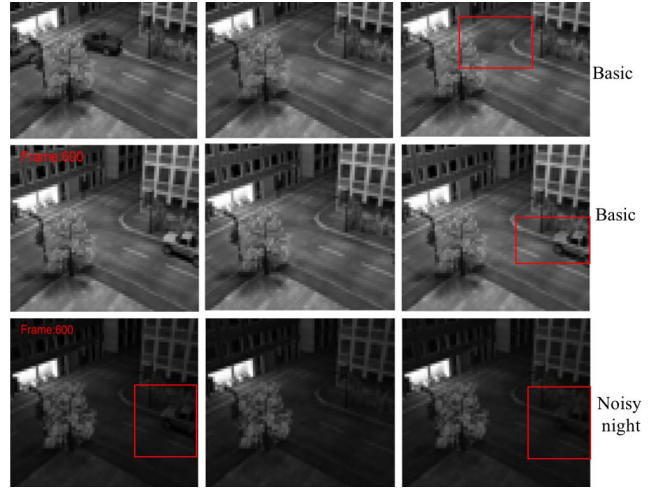}
 	\caption{Qualitative analysis of the background estimated by sWLR and APG on the~{\it Basic} and~{\it Noisy night} scenario.~Frame number 600 has static foreground in both scenarios.~APG can not remove the static foreground object from the background.~On the other hand, in frame number 210, the low-rank background estimated by APG has still some black patches.~In both scenarios sWLR can completely remove the static foreground.}\label{basic}
 \end{figure}

This special choice of the weight is justified as follows:~in background subtraction, we only need to put large weights on the columns (frames) that are more likely to be the background and leave the rest of the columns unweighted~(and thus with weight $1$). Our new algorithm is not based on matrix factorization to address the rank constraint~\cite{duttali_cvpr}. Instead, we exploit the dependence of ${X}_2$ on ${X}_1$ in the optimal solution.
We will use~Theorem~\ref{theorem 7} to device an iterative process to solve~(\ref{2nd approximate golub's problem}).
We assume that ${\rm r}({X}_1)=k.$ Then any ${X}_2$ such that ${\rm r}({X}_1\;\;{X}_2)\le r$ can be given in the form
\vspace{-0.1in}
$${X}_2={X}_1C+D,$$ 
~\\[-0.21in]
for some arbitrary matrices $C\in\mathbb{R}^{k\times (n-k)}$ and $D\in\mathbb{R}^{m\times (n-k)}$, such that ${\rm r}(D)\le r-k$. Therefore,~(\ref{2nd approximate golub's problem}) becomes an constrained weighted low-rank approximation problem:
\begin{align}\label{main problem 2}
\min_{\substack{{X}_1,C,D\\{\rm r}(D)\le r-k}}\left(\|(A_1-{X}_1)\odot W_1\|_F^2+\|A_2-{X}_1C-D\|_F^2\right).
\end{align}
~\\[-0.12in]
Denote $F({X}_1,C,D)=\|(A_1-{X}_1)\odot W_1\|_F^2+\|A_2-{X}_1C-D\|_F^2$ as the objective function.
Assume that at the $p$-th step we have $(X_1)_p$. We need to find $(C_p,D_p)$ by solving
$$
\min_{C,D} F((X_1)_p,C,D).
$$
Then Theorem~\ref{theorem 7} suggests
~\\[-0.12in]
$$
(X_1)_pC_p=P_{(X_1)_p}(A_2)~{\rm and}~D_p=H_{r-k}(P_{(X_1)_p}^\perp(A_2)).$$
So, if $(X_1)_p$ has its QR decomposition:
~\\[-0.12in]
$$
 (X_1)_p=Q_pR_p,
 $$
 then
 $$
 C_p=R_p^{-1}Q_p^TA_2
 $$
 and
~\\[-0.05in]
 $$
 D_p=H_{r-k}((I-Q_pQ_p^T)A_2)=U_p(\Sigma_p)_{r-k}V_p^T,
 $$
with $U_p\Sigma_pV_p^T$ being a SVD of $(I-Q_pQ_p^T)A_2$.
 
We are only left to find $({X}_1)_{p+1}$ given $(C_p,D_p)$ via the following iterative scheme:
\begin{align}\label{update rule}
\displaystyle{({X}_1)_{p+1}=\arg\min_{{X}_1}F({X}_1,C_p,D_p)}.
\end{align}
We will update ${X}_1$ row-wise. Let ${X}_1(i,:)$ denote the $i$-th row of the matrix ${X}_1$. We set $\frac{\partial}{\partial{X}_1}F({X}_1,C_{p},B_{p})=0$~and obtain
$$
-(A_1-{X}_1)\odot W_1\odot W_1-(A_2-{X}_1C_{p}-D_{p})C_{p}^T=0.
$$
Solving the above expression for ${X}_1$ sequentially along each row produces
\begin{align*}
({X}_1(i,:))_{p+1}=&(E(i,:))_p({\rm diag}(W_1^2(i,1)&\\&\;\;W_1^2(i,2)\cdots W_1^2(i,k))+C_{p}C_{p}^T)^{-1},
\end{align*}
where $E_p=A_1\odot W_1\odot W_1+(A_2-D_{p})C_{p}^T$. Therefore, we have the following algorithm.
\begin{algorithm}\label{algorithm1}
	\SetAlgoLined
	\SetKwInOut{Input}{Input}
	\SetKwInOut{Output}{Output}
	\SetKwInOut{Init}{Initialize}
	\nl\Input{$A=(A_1\;\;A_2)\in\mathbb{R}^{m\times n}$ (the given matrix); $W=(W_1\;\;\mathbf{1})\in\mathbb{R}^{m\times n},$ (the weight); threshold  $\epsilon>0$\;}
	\nl\Init {$(X_1)_0$\;}
	\BlankLine
	\nl \While{not converged}
	{
		\nl $({X}_1)_p=Q_pR_p$, $(I_m-Q_pQ_p^T)A_2=U_p\Sigma_pV_p^T$\;
		\nl $C_{p}=R_p^{-1}Q_p^TA_2 $\;
		\nl $D_{p}=U_p(\Sigma_p)_{r-k}V_p^T$\;
		\nl $E_p=A_1\odot W_1\odot W_1+(A_2-D_{p})C_{p}^T$\;
			\nl \For {$i=1:m$}
			{
		\nl $({X}_1(i,:))_{p+1}=(E(i,:))_p({\rm diag}(W_1^2(i,1)
		\myindent{0.3} \;\;W_1^2(i,2)\cdots W_1^2(i,k))+C_{p}C_{p}^T)^{-1}$\;
	}
		\nl $p=p+1$\;
	}
	\BlankLine
	\nl \Output{$(X_1)_{p}, (X_1)_{p}C_{p}+D_{p}.$}
	\caption{sWLR Algorithm}
\end{algorithm}

The update rule for Algorithm~sWLR is \vspace{-0.1in}$${X}_{p+1}=(({X}_1)_{p}\;\; ({X}_1)_{p}C_{p}+D_{p}),$$ 
~\\[-0.15in]
with ${\rm r}(({X}_1)_{p})=k$, ${\rm r}(({X}_1)_{p}C_p)\leq k$, ${\rm r}(D_{p})\leq r-k$, and so, ${\rm r}({X}_{p+1})\leq r$.
\vspace{-0.05in}
\section{Numerical Experiments}\label{numerical}
\vspace{-0.05in}
In this section we will demonstrate the performance of our algorithm in solving the background estimation problem and compare it with WLR in~\cite{duttali_cvpr} and the RPCA methods~\cite{LinChenMa,APG,candeslimawright}.

\vspace{-0.01in}
\subsection{Implementation Details}
\vspace{-0.01in}
Let $X_{sWLR}=({X}_1^*\;\; {X}_1^*C^*+ D^*)$ where $({X}_1^*, C^*$, $D^*)$ is a solution to~(\ref{main problem 2}).~We denote $X_p$ as our approximation to $X_{sWLR}$ at $p$th iteration. Recall that $X_p=(({X}_1)_{p}\;\; ({X}_1)_{p}C_p+D_p).$ We denote $\|X_{p+1}-X_{p}\|_F=Error_p$ and use $\frac{\|Error_p\|_F}{\|X_{p}\|_F}$ as a measure of the relative error. For a threshold $\epsilon>0$ the stopping criteria of our algorithm at the $(p+1)$th iteration is  $\|Error_p\|_F<\epsilon$ or $\frac{\|Error_p\|_F}{\|X_{p}\|_F}<\epsilon$ or if the maximum iterations attained. The algorithm performs the best when we initialize ${X}_1$ as a random matrix and takes 5--10 iterations to converge.

Recall that, the {\it Robust PCA}~(RPCA) method for background estimation problems uses the fact that the background frames, $X$, have a low-rank structure and the foreground $A-X$ is sparse~\cite{LinChenMa,APG,candeslimawright} and solves:
\vspace{-0.08in}
\begin{equation}\label{rpca}
\min_X\{\|A-X\|_{\ell_1}+\lambda \|X\|_*\}.
\vspace{-0.04in}
\end{equation}
\noindent For RPCA, we use the inexact augmented Lagrange multiplier~(iEALM) method proposed by Lin et. al.\ ~\cite{LinChenMa}, and the accelerated proximal gradient~(APG) algorithm proposed by Wright et. al.\ ~\cite{APG}.~For iEALM and APG we set $\lambda={1}/{\sqrt{{\rm max}\{m,n\}}}$, and for iEALM we choose $\mu=1.5, \rho=1.25$~\cite{LinChenMa,candeslimawright,APG}.~A threshold equal to $10^{-7}$  is set for all algorithms.

\vspace{-0.01in}
\subsection{Experimental Setup}
\vspace{-0.1in}
We perform our experiments on~the Stuttgart synthetic video data set~\cite{cvpr11brutzer}.~It is a computer generated video sequence, that comprises gradual or sudden change of illumination, a dynamic background containing non-stationary objects and a static foreground, camouflage, and sensor noise or compression artifacts.~We perform qualitative and quantitative analysis on two different test scenarios of the sequence:
(i) {\it Basic} and~(ii) {\it Noisy night}. Each scenario has 600 frames with identical foreground and background objects. Frame numbers 551 to 600 have static foreground, and frame numbers 6 to 12 and 483 to 528 have no foreground.~Additionally, the foreground comes with high quality ground truth mask for each video frame.
\subsection{Comparison between RPCA and sWLR}
\vspace{-0.05in}
In Figure~\ref{weight_435}, we demonstrate the effects of using big weights on the frames in the first block $X_1$ for sWLR. With proper choice of $r$ and $k$ the large weights in $W_1$ produces a better background estimation.~In Figure~\ref{basic}, we present qualitative comparison between the background estimated by sWLR and RPCA algorithms. Since APG and iEALM both have same reconstruction we only present APG here.~In both scenarios, sWLR provides a substantially better background estimation than APG.

For quantitative comparison between different methods we use~the most advanced measure structural similarity index~(SSIM) in Figure~\ref{ssim}(a) and~(b).~According to~\cite{mssim},~the SSIM index can be viewed as a robust quality measure of a recovered image, compare to the other image that is regarded as of perfect quality.~It compares the luminance change, contrast change, and structural change in the recovered image and agree with human visual perception the most compare to any other standard measures.~Our background estimation experiments for sWLR are based on a prior knowledge of the background frame indexes.~From the ground truth, we know that the entire sequence has 60 foreground frames that has less than 10 pixels.~Given 60 pure background frames we choose $k=\Bigl\lceil 60/i_1\Bigr\rceil$  by random sampling, where $i_1\in\{1,2,3,4\}$.~We~set $r=k+1$. For {\it Basic} and {\it Noisy night} $k=15$ and $k=20$ respectively, are the best choice.~Since the qualitative and quantitative results for background estimation are same for both WLR and sWLR we only provide their runtime comparison~(see Figure~\ref{ssim}~(c)).~In Figure~\ref{ssim}(c), the increment in time for the RPCA algorithms as we increase the dimension of the test matrix by adding more frames can be attributed by computation of a larger rank SVD in each step of their iteration. On the other hand, sWLR performs a fixed rank SVD once the background frames are learned. 

\section{Conclusion}
\vspace{-0.01in}
In this paper we presented a simple and fast numerical algorithm to solve a weighted low rank approximation problem for a special family of weights.~To demonstrate its use in the real world problems we performed background estimation from video sequences when a prior knowledge of approximated background fames is available.~We did not address the question on how to automatically learn the weight from the data which is treated in~\cite{duttali_cvpr}.~The performance of our weighted low-rank approximation algorithm over the existing RPCA algorithms shows the fact that a weighted Frobenius norm can be made robust to sparse outliers.~With additional knowledge of some approximate frames learned from the data, our algorithm can outperform the RPCA algorithms in terms of accuracy and efficiency.


\end{document}